\documentclass{article}
\usepackage[left=3cm,right=3cm,top=3cm,bottom=2cm]{geometry} 
\usepackage{amsmath} 
\usepackage{amssymb}
\usepackage{cite}
\usepackage{color}

\usepackage[T1,T2A]{fontenc}
\usepackage[utf8]{inputenc}

\begin{document}

\begin{Large}
\centerline{On a problem related to "second" best approximations to a real number}
\vskip+0.5cm
\centerline{Pavel Semenyuk}
\end{Large}

\vskip+1cm
\section{Introduction}

\vskip+0.3cm
For a given irrational number $\alpha$ we consider its continued fraction expansion 
\begin{equation}\label{q}
	\alpha = [a_0, a_1, a_2, ...] = a_0 + \dfrac{1}{a_1+\dfrac{1}{a_2+...}}, ~a_0\in\mathbb{Z}, ~a_i\in\mathbb{Z}_{+}, ~i = 1, 2, ...
\end{equation}
 and its convergents $\dfrac{p_n}{q_n} = [a_0, a_1, ..., a_n]$. We call two irrational numbers  $\alpha = [a_0, a_1, a_2, ...]$ and  $\beta = [b_0, b_1, b_2, ...]$ equivalent  and write $\alpha\sim\beta$, if 
there exist positive integers $m$ and $n$ such that 
$a_{n+k} = b_{m+k}, k=1,2,3,...$, that is the tails of continued fraction expansions for $\alpha$ and $\beta$ coincide.
Let
$$
\psi_{\alpha}(t) = \min\limits_{1\leqslant q \leqslant t, q\in\mathbb{Z}} ||q\alpha|| ,\,\,\,\,\,\text{where}\,\,\,\,\,  ||x|| = \min\limits_{n\in\mathbb{Z}}|x-n|
$$ 
be the irrationality measure function for $\alpha$.
By Lagrange's theorem on best approximations (see \cite{khinchin1964}) it is  a piecewise constant function, namely
$$ \psi_{\alpha}(t) = |q_{\nu}\alpha - p_{\nu}|
\,\,\,\,\,\text{
for}
\,\,\,\,\,
q_\nu\leqslant t\leqslant q_{\nu+1}.
$$
 Many classical result concerning rational approximations to a real number $\alpha$ can be formulated in terms of $\psi_{\alpha}(t)$. 
 For example, in terms of irrationality measure function $\psi_\alpha$ one can define the Lagrange constant
 \begin{equation}\label{l}
 \lambda(\alpha) = 
 \liminf_{t\to \infty} t\cdot \psi_\alpha (t).
 \end{equation}
 The set of all possible values of $\lambda(\alpha)$ is known as 
 Lagrange spectrum
$$
\mathbb{L} = \{\lambda ~|~\exists\alpha\in\mathbb{R\setminus Q} \colon \lambda = \lambda(\alpha)  \}.
$$
It is very well studied 
(see for example  \cite{cusick1989markoff}).

\vskip+0.3cm

In \cite{moshchevitin2017} Moshchevitin considered an irrationality measure function
$$
	\psi_{\alpha}^{[2]}(t) = \min\limits_{\substack{(q, p)\colon q, p \in\mathbb{Z}, 1\leqslant q\leqslant t, \\ (p, q) \neq (p_n, q_n) ~\forall n\in\mathbb{Z_{+}}}} |q\alpha -p|,
$$
related to so-called "second best" approximations,
the corresponding Diophantine constant
\begin{equation}\label{k}
 \mathfrak{k}(\alpha) = \liminf\limits_{t\to\infty} t \cdot\psi_{\alpha}^{[2]}(t) 
\end{equation}
and the corresponding spectrum 
$$
\mathbb{L}_2 = \{ \lambda ~| ~\exists\alpha\in\mathbb{R\setminus Q} \colon \lambda = \mathfrak{k}(\alpha)  \}.
$$ 
In particular, he 
studied different properties of the function $\psi_{\alpha}^{[2]}(t)$
and
calculated two largest elements of the spectrum $\mathbb{L}_2$.  In the present paper we calculate the value for the third element of the spectrum $\mathbb{L}_2$.

\vskip+0.3cm

\section{Results on the spectrum $\mathbb{L}_2$}

\vskip+0.3cm

The following structural result about $\mathbb{L}_2$ was proven
in \cite{moshchevitin2017}. 

\vskip+0.3cm

\textbf{Theorem 1.} 
\begin{itemize}

\item[1.] The largest element of $\mathbb{L}_2$ is $\dfrac{4}{\sqrt{5}}$. Moreover, $\mathfrak{k}(\alpha) = \dfrac{4}{\sqrt{5}}$ if and only if $\alpha\sim\dfrac{1+\sqrt{5}}{2} = [1; \overline{1}]$;
\item[2.] If $\alpha\not\sim\dfrac{1+\sqrt{5}}{2}$, then $\mathfrak{k}(\alpha)\leqslant\dfrac{4}{\sqrt{17}}$;
\item[3.] The equality $\mathfrak{k}(\alpha) = \dfrac{4}{\sqrt{17}}$ holds  if and only if $\alpha\sim\dfrac{1+\sqrt{17}}{2} = [2; \overline{1,1,3}]$;
\item[4.] The point $\dfrac{4}{\sqrt{17}}$ is an isolated point of $\mathbb{L}_2$.
\item[5.] The whole segment $\left[ 0, \frac{12}{21+\sqrt{15}}\right]$ belongs to $\mathbb{L}_2$.

\end{itemize}
\vskip+0.3cm

The main result of the present paper is as follows.

\vskip+0.3cm

\textbf{Theorem 2.} 
\begin{enumerate}
\item If $\alpha$ is irrational and not equivalent neither to  $\frac{1+\sqrt{5}}{2}$ nor to $\frac{1+\sqrt{17}}{2}$, then $\mathfrak{k}(\alpha) \leqslant \frac{164}{13\sqrt{173}}$.
\item Let $\alpha_0 = \frac{13\sqrt{173} + 39}{82} = [2; \overline{1, 1, 3, 1, 1, 1, 1, 3}]$.
Then $\mathfrak{k}(\alpha) =\mathfrak{k}(\alpha_0) = \frac{164}{13\sqrt{173}}$ if and only if $\alpha \sim \alpha_0$.
\item The point $\frac{164}{13\sqrt{173}}$ is an isolated point of $\mathbb{L}_2$.
\end{enumerate}

\vskip+0.3cm

\section{Auxiliary results and notation}

\vskip+0.3cm

Consider the tails 
 $$
 \alpha_n = [a_n; a_{n+1}, a_{n+2}, a_{n+3}, ...],\,\,\,\,\,
 \alpha^{*}_n = [0; a_n, a_{n-1}, a_{n-2}, ..., a_1]
 $$
 of continued fraction  (\ref{q}). In terms of these quantities it is natural to give a formula 
 \begin{equation}\label{pe}
 \left| \alpha - \frac{p_n}{q_n}\right| = \frac{1}{q_n^2 (\alpha_{n+1}+\alpha_n^*)}
 \end{equation}
 for the approximation to $\alpha$ by its convergent fraction.
 Relation (\ref{pe}) is known as Perron's formula. By means of (\ref{pe}) one can express the value of $\lambda(\alpha)$ 
 defined in (\ref{l}) as
  \begin{equation}\label{pe1}
  \lambda(\alpha) = \liminf_{n\to \infty} \frac{1}{{\alpha_{n+1}+\alpha_n^*}}.
 \end{equation}
 An analog of the formula (\ref{pe1}) for the value of $\frak{k}(\alpha)$ was obtained in  \cite{moshchevitin2017} .
 It is as follows.
 Consider  the  quantities
$$
\varkappa^1_n(\alpha) = \dfrac{(1 + \alpha^{*}_{n-1})(\alpha_n - 1)}{\alpha_n + \alpha^{*}_{n-1}},\,\,\,\,\,\,
\varkappa^2_n(\alpha) = \dfrac{(1 - \alpha^{*}_n)(\alpha_{n+1} + 1)}{\alpha_{n+1} + \alpha^{*}_n}
,
\,\,\,\,\,\,
\varkappa^4_n(\alpha) = \dfrac{4}{\alpha_n + \alpha^{*}_{n-1}}.
$$
(we follow the notation from  \cite{moshchevitin2017}).
Then
 if $\alpha \not\sim \frac{1+\sqrt{5}}{2}$ one has 
\begin{equation}\label{one}
\mathfrak{k}(\alpha) = \liminf\limits_{n\to\infty\colon a_n\geqslant 2}\min(\varkappa^1_n, \varkappa^2_n, \varkappa^4_n)
\end{equation}
(see Hilfssatz 13 from  \cite{moshchevitin2017}).
We should note that it is clear that for $\alpha \sim \beta$ one has $ \mathfrak{k}(\alpha) = \mathfrak{k}(\beta) $.

\vskip+0.3cm

Besides formula (\ref{one}) we need some other auxiliary statements from \cite{moshchevitin2017}.
We formulate them as the following

\vskip+0.3cm
 
\textbf{Proposition 1.} Let $\alpha$ be an irrational number not equivalent to $\dfrac{1+\sqrt{5}}{2} = [1; \overline{1}]$ and $\dfrac{1+\sqrt{17}}{2} = [2; \overline{1,1,3}]$.
\begin{itemize}

\item[1.] If for infinitely many $n$ in continued fraction expansion (\ref{q})  one has $a_n \geqslant 5$, then $\mathfrak{k}(\alpha) \leqslant \frac{4}{5}$;
\item[2.] If for all $n$ sufficiently large   in continued fraction expansion (\ref{q})  one has $a_n \leqslant 4$ and for infinitely many $n$  the equality $a_n = 4$ occurs, then $\mathfrak{k}(\alpha) \leqslant \frac{4}{3 + \sqrt{2}}$;
\item[3.] If for all $n$ sufficiently large $a_n \leqslant 4$ and for infinitely many $n$  the equality $a_n = 2$ occurs, then $\mathfrak{k}(\alpha) \leqslant \sqrt{2} - \frac{1}{2}$;
\item[4.] If for infinitely many $n$  one has $a_n \in \{1, 3\}  $ and
\item[4.1.]  for infinitely many $n$  one has $a_{n-1} = a_n = 3$, then $\mathfrak{k}(\alpha) \leqslant \frac{39}{43}$;
\item[4.2.] for infinitely many $n$  one has $a_{n-1} = 3, a_n = 1, a_{n+1} = 3$, then $\mathfrak{k}(\alpha) \leqslant \frac{136}{145}$;
\item[4.3.]  for infinitely many $n$ one has $a_{n-2} = a_{n-1} = 1, a_n = 3, a_{n+1} = a_{n+2} = a_{n+3} = 1$, then $\mathfrak{k}(\alpha) \leqslant \frac{180}{187}$  
\end{itemize}

\vskip+0.3cm

{\bf Remark.} Although the main result on the structure of $\mathbb{L}_2$ from the paper 
\cite{moshchevitin2017} 
is correct, its  proof  from 
\cite{moshchevitin2017} contains 
 a mistake:
in the cases 4.2 and 4.3 in \cite{moshchevitin2017} 
instead of upper bounds $\frac{136}{145}$ and $ \frac{180}{187}$  correspondingly,   different incorrect bounds were given.\\
\vskip+0.3cm

The bounds for $\mathfrak{k}(\alpha)$  from Proposition 1 can be written as the following table.
\vskip+0.3cm

\begin{center}
\begin{tabular}{ | c | c | c | c | }
\hline
Case & Conditions for partial quotients & Upper bound & Numerical value\\ \hline
1 & $\text{inf. many}\,\, a_n \geqslant 5$ & $4/5$ & 0.8 \\ \hline
2 & $\text{almost all}\,\, a_n\le 4, \,\,\text{inf. many}\,\, a_n = 4$ & $4/(3+\sqrt{2})$ & $0.906163^+$ \\ \hline
3 & $\text{almost all}\,\, a_n\le 4, \,\,\text{inf. many}\,\, a_n = 2$  & $\sqrt{2} - 1/2$ & $0.914213^+$ \\ \hline
4.1 & $\text{inf. many patterns}\,\,$  3 3 & 39/43 & $0.906976^+$ \\ \hline
4.2 &$\text{inf. many patterns}\,\,$   3 1 3 & 136/145 & $0.937931^+$ \\ \hline
4.3 & $\text{inf. many patterns}\,\,$   1 1 3 1 1 1 & 180/187 & $0.962566^+$ \\
\hline
\end{tabular}
\end{center}

\vskip+0.3cm

\section{The value of $\mathfrak{k} (\alpha_0)$}
  
\vskip+0.3cm

In this section we prove the following
  
\vskip+0.3cm

\textbf{Lemma 1.}  $\mathfrak{k}(\alpha_0) = \liminf\limits_{n\to\infty\colon a_n\geqslant 2}\varkappa^4_n = \dfrac{164}{13\sqrt{173}}$.

\vskip+0.3cm

\textbf{Proof.} Consider blocks of digits $ {A} = \{1,1,3\}$ and ${ B} = \{1,1,1,1,3\}$. 
Now $\alpha_0 = [2; \overline{1, 1, 3, 1, 1, 1, 1, 3}]$ can be written as
$\alpha_0 =[2;\overline{A B}]$.
In view of (\ref{q}), it suffices to prove that for all large enough $n$ such that $a_n = 3$,  the value of $\varkappa^4_n$ is less than both $\varkappa^1_n$ and $\varkappa^2_n$. Let us consider separately those values of  $n$ which are either $\equiv 0 \pmod{8}$ (if we take a 3 after a block of four 1's) or  $\equiv 3 \pmod{8}$ (if we take a 3 after a block of two 1's). Thus, we need to consider two possible limits for each $\varkappa^1_n, \varkappa^2_n, \varkappa^4_n$.

Let us calculate all six of them:
$$
\lim\limits_{k\to\infty} \varkappa^1_{8k-5} = \dfrac{(1+[0;\overline{AB}])([3;\overline{BA}]-1)}{[0;\overline{AB}]+[3;\overline{BA}]} = \dfrac{167}{13\sqrt{173}} = 0.976674^+,
$$
$$
\lim\limits_{k\to\infty} \varkappa^1_{8k} = \dfrac{(1+[0;\overline{BA}])([3;\overline{AB}]-1)}{[0;\overline{BA}] + [3;\overline{AB}]} = \dfrac{169}{13\sqrt{173}} = 0.988371^+,
$$
$$
\lim\limits_{k\to\infty} \varkappa^2_{8k-5} = \dfrac{(1-[0;3,\overline{AB}])([1;1,1,1,3,\overline{AB}]+1)}{[0;3,\overline{AB}]+[1;1,1,1,3,\overline{AB}]} = \dfrac{169}{13\sqrt{173}} = 0.988371^+,
$$
$$
\lim\limits_{k\to\infty} \varkappa^2_{8k} = \dfrac{(1-[0;3,\overline{BA}])([1;1,3,\overline{BA}]+1)}{[0;3,\overline{BA}]+[1;1,3,\overline{BA}]} = \dfrac{167}{13\sqrt{173}} = 0.976674^+,
$$
$$
\lim\limits_{k\to\infty} \varkappa^4_{8k-5} = \lim\limits_{k\to\infty} \varkappa^4_{8k} = \dfrac{4}{[0;\overline{AB}]+[0;\overline{BA}]+3} = \dfrac{164}{13\sqrt{173}} = 0.959129^+.
$$
We see that 
the last one 
is the smallest one. This proves Lemma 1.

 \vskip+0.3cm

  \section{Extremality of $\mathfrak{k} (\alpha_0)$}
  
  \vskip+0.3cm

  In this section we finalise the proof of Theorem 2.
    \vskip+0.3cm

In cases 1-3, 4.1 and 4.2 of Proposition 1  we have $\mathfrak{k}(\alpha) \leqslant \frac{136}{145} < \mathfrak{k}(\alpha_0) = \frac{164}{13\sqrt{173}}= 0.959129^+$,
meanwhile the bound of the case 4.3 satisfies the inequality $ \frac{180}{187} >\mathfrak{k}(\alpha_0)$.
So we should consider  in more details just the case 4.3. First of all we consider the following four subcases in the case 4.3:

  \vskip+0.3cm

\begin{itemize}
\item[4.3.1] One has infinitely many patterns $\{1, 1, \boldsymbol{3}, 1, 1, 1, 1, 1\}$ (here and in other cases the bold \textbf{3} represents the $n$th element of the continued fraction). In that case
\begin{multline*}
\mathfrak{k}(\alpha) \leqslant \liminf\limits_{n\to\infty\colon a_n\geqslant 2}\varkappa^4_n(\alpha) \leqslant \dfrac{4}{[0;1,1,\overline{3,1}] + [3;1,1,1,1,1,\overline{1,3}]} =\\= \dfrac{4}{\frac{\sqrt{21}+1}{10} + \frac{\sqrt{21}+943}{262}} = 0.958087^+.
\end{multline*}
\item[4.3.2] One has infinitely many patterns $\{1, 1, \boldsymbol{3}, 1, 1, 1, 3\}$. In that case
\begin{multline*}
\mathfrak{k}(\alpha) \leqslant \liminf\limits_{n\to\infty\colon a_n\geqslant 2}\varkappa^4_n(\alpha) \leqslant \dfrac{4}{[0;1,1,\overline{3,1}] + [3;1,1,1,3,1,1,\overline{3,1}]} =\\= \dfrac{4}{\frac{\sqrt{21}+1}{10} + \frac{\sqrt{21}+4575}{1258}} = 0.952692^+.
\end{multline*}
\item[4.3.3] One has infinitely many patterns $\{ABB\}$ (or, in other words, infinitely many patterns \\$\{1, 1, 3, 1, 1, 1, 1, \boldsymbol{3}, 1, 1, 1, 1, 3, 1, 1\}$. In that case
\begin{multline*}
\mathfrak{k}(\alpha) \leqslant \liminf\limits_{n\to\infty\colon a_n\geqslant 2}\varkappa^4_n(\alpha) \leqslant \dfrac{4}{[0;1,1,1,1,3,1,1,\overline{1,3}]+[3;1,1,1,1,3,1,1,\overline{1,3}]} =\\= \dfrac{4}{\frac{3329+\sqrt{21}}{5470} + \frac{19739+\sqrt{21}}{5470}} = 0.948123^+.
\end{multline*}
\item[4.3.4] One has infinitely many patterns $\{AAB\}$ (or, in other words, infinitely many patterns \\$\{1, 1, 3, 1, 1, 1, 1, \boldsymbol{3}, 1, 1, 3, 1, 1, 3, 1, 1\}$. In that case
\begin{multline*}
\mathfrak{k}(\alpha) \leqslant \liminf\limits_{n\to\infty\colon a_n\geqslant 2}\varkappa^4_n(\alpha) \leqslant \dfrac{4}{[0;1,1,1,1,3,1,1,\overline{1,3}]+[3;1,1,3,1,1,3,1,1,\overline{3,1}]} =\\= \dfrac{4}{\frac{3329+\sqrt{21}}{5470} + \frac{120383+\sqrt{21}}{33802}}= 0.959006^+.
\end{multline*}

\end{itemize}

We see that in all the subcases  4.3.1, 4.3.2, 4.3.3, 4.3.4 the upper bound for $\mathfrak{k}(\alpha) $ obtained is strictly less than $\mathfrak{k}(\alpha_0) = \frac{164}{13\sqrt{173}}$.

\medskip

Let us describe all the cases considered. 
Cases 1-3 cover all the continued fractions with infinitely many terms different from  1 and  3. Cases 4.1, 4.2, 4.3.1 and 4.3.2 cover all the fractions with infinitely many patterns of the form $[..., 3, \underbrace{1, ..., 1}_{k}, 3, ...]$, where $k\not\in\{2, 4\}$.
Therefore, the only fractions which we did not consider in the cases above, eventually consist of blocks of the form $ {A} $ and ${ B}$. Case 4.3.3 covers all the fractions with infinitely many patterns $ABB$, and case 4.3.4 covers fractions with infinitely many patterns $AAB$. Thus, the only fractions left in consideration are eventually  periodic with the period $\overline{AB}$, or, in other words, are equivalent to $\alpha_0=\frac{13\sqrt{173}+39}{82}$. Together with Lemma 1 this  proves Theorem 2.

\end{document}